\documentclass[12pt,leqno]{article}
\usepackage{latexsym,amsmath,amssymb,epsfig}
\title{The Classification of Rank $4$ Locally Projective Polytopes and Their Quotients
{ }}
\author{Michael I.~Hartley\\ {\small The University of Nottingham}\\
{\small Malaysian Campus}
{\small 2 Jalan Conlay}\\
{\small Kuala Lumpur, 50450, Malaysia}\\
{\small Michael.Hartley@nottingham.edu.my}}
\date{}

\advance\voffset by-20mm
\makeatletter
\def\@eqnnum{\hbox to .01pt{}\rlap{\bf \hskip -
\displaywidth\theequation}}
\makeatother
\renewcommand{\theequation}{\thesection.\arabic{equation}}
\newcommand{\seq}{\setcounter{equation}{0}}
\newenvironment{proof}{\par \noindent {\em
Proof\/}:\quad}{\hspace*{\fill} $\Box$ \par \vspace*{1ex}}
\newcommand{\bgpf}{\begin{proof}}
\newcommand{\ndpf}{\end{proof}}
\newenvironment{reason}{\par \noindent {\em
Reason\/}:\quad}{\hspace*{\fill} $\Box$ \par \vspace*{1ex}}
\newcommand{\bgrs}{\begin{reason}}
\newcommand{\ndrs}{\end{reason}}

\newtheorem{theorem}[equation]{\hspace{-0.3em}}
\newcommand{\bgth}{\begin{theorem}}
\newcommand{\ndth}{\end{theorem}}

\newcommand{\thm}{{\sc Theorem}\quad}
\newcommand{\cor}{{\sc Corollary}\quad}

\newcommand{\bgeq}{\begin{equation}}
\newcommand{\bgqy}{\begin{eqnarray*}}
\newcommand{\bgry}{\begin{array}}
\newcommand{\ndeq}{\end{equation}}
\newcommand{\ndqy}{\end{eqnarray*}}
\newcommand{\ndry}{\end{array}}
\def\Aut{\mathop{{\rm Aut}}}

\newcommand{\CH}{{\cal H}}
\newcommand{\CI}{{\cal I}}
\newcommand{\CK}{{\cal K}}
\newcommand{\CL}{{\cal L}}

\newcommand{\CP}{{\cal P}}
\newcommand{\CQ}{{\cal Q}}

\newcommand{\scl}[1]{\langle{#1}\rangle}
\begin{document}
\maketitle
\centerline{\it In Memory of H. S. M. ``Donald'' Coxeter, 1907--2003}
\begin{abstract}
This article announces the completion of the classification of
rank 4 locally projective polytopes and their quotients.
There are seventeen universal locally projective polytopes (nine nondegenerate).
Amongst their 441 quotients are a further four (nonuniversal)
regular polytopes, and 152 nonregular but section regular polytopes.
All 156 of the latter have hemidodecahedral facets or hemi-icosahedral
vertex figures. It is noted that, remarkably, every rank 4 locally projective
section regular polytope is finite.
The article gives a survey of the literature of locally projective polytopes
and their quotients, and fills one small gap in the classification in rank 4.
\end{abstract}

%
%
%This is file 1.tex (version 15.07.1996).
%
%
\section{Introduction}
\label{intro}
\seq

The classical study of polytopes concentrated on what are now called ``spherical'' polytopes.
The modern ``abstract'' polytopes include tesselations of euclidian and hyperbolic space and
other spaceforms, as well as objects are best described by their ``local'' topology, that is,
the topology of their facets and vertex figures, or smaller sections.
Good accounts of the classical theory and of the emergence of the theory of
the abstract polytopes are to be found in \cite{cox} or \cite{msbook}.

Extensive study has been
made by McMullen, Schulte and others of the so called ``locally toroidal'' polytopes, that is,
polytopes whose minimal nonspherical sections are toroids. See for example Chapters 10 to 12
of \cite{msbook}, which bring together and update earlier work on the topic
(such as \cite{mstwis}).
Comparatively little has been done, until recently, on polytopes of other local topological types.
The next three paragraphs summarise briefly the literature on locally projective polytopes.

There are two non-equivalent definitions used in the literature for a ``locally $X$'' polytope.
The broader definition is that the minimal nonspherical sections of the polytope must
be of topological type $X$ (see for example \cite{mslocproj}).
The narrower definition is that these minimal nonspherical sections
must actually be the facets or vertex figures (or both) of the polytope
(see for example \cite{mstwis}). For rank $4$ polytopes,
both definitions are equivalent: there do not exist
nonspherical polytopes of rank less than $3$,
so the minimal nonspherical sections of a rank $4$ polytope must be
of rank $3$, if they exist at all. The terminology ``locally X''
is normally restricted to regular polytopes, or at least to section
regular polytopes. (A {\it section regular} polytope is one for which
the isomorphism type of its sections $F/G$ depends only on the ranks of
$F$ and $G$.)

In the late 1970's and early 1980's, Gr\"unbaum and Coxeter independently discovered a
self-dual polytope with 11 hemi-icosahedral facets
(see \cite{cox11} and \cite{gru11}). During the same period, (see \cite{cox57}),
Coxeter discovered another
self-dual polytope with 57 hemidodecahedral facets. 

As the theory of abstract polytopes developed, these sporadic examples were joined
by others that arose out of more general considerations. 
Schulte, in \cite{amalg}, noted some examples
of locally projective polytopes of type $\{4,3,4\}$, that is, with hemicubes for facets, 
and octahedra
or hemioctahedra for vertex figures. Some years later, in 1994,
McMullen in \cite{mslocproj} examined 
locally projective polytopes of type $\{3^{k-2},4,3,3,4,3^{n-k-3}\}$ (that is, $\{3,4,3^{n-3}\}$,
$\{3,3,4,3^{n-4}\}$, $\{4,3,3,4,3^{n-5}\}$ respectively for $k=0,1,2$), where $n\geq5$.
This was the first paper to give an infinite family of locally projective polytopes.
As was shown in \cite{har5}, the theory developed by Schulte in \cite{amalg}
can be easily applied to produce another infinite family. However, Schulte had chosen a
different focus for that paper.

The most recent work on the problem has been a series of articles by the author and others 
(see \cite{har4}, \cite{har5} and \cite{harlee}) which focused on specific examples
of locally projective polytopes and their quotients. That work
almost completed the classification of 
locally projective rank $4$ polytopes and their quotients. This article serves
mainly as a survey of the classification (see Section~\ref{class}), but also (in
Section~\ref{uncases}) closes the one small gap that remains.
Section~\ref{formscope} contains observations that show that
the cases examined are indeed all that can occur.

The most
remarkable feature of the classification is the fact that all such polytopes are finite. This
is in contrast with the rank 5 case, where infinite locally projective polytopes exist (an example 
being $\{\{4,3,3\},\{3,3,5\}/2\}$ (see Theorem 5.2 of \cite{har5}).

\section{The Form and Scope of the Classification.}
\label{formscope}
\seq

The usual way to search for polytopes with given (regular) facets
$\CK$ and vertex figures $\CL$ is to find the {\it universal} such
polytope $\CP=\{\CL,\CK\}$, and then seek its quotients. The universal
polytope is characterised by its automorphism group
$\Aut(\CP)=\Gamma(\CP)=[\CL,\CK]$. By definition, all
polytopes of the desired form will be quotients of $\CP$.
It is guaranteed that if any polytope exists with the desired
facets and vertex figures, then the universal polytope $\CP$ does
(see Theorem 4A2 of \cite{msbook} or Theorem 2.5 of \cite{har6}).

A locally projective section regular polytope of rank $4$ must have spherical or projective facets
and vertex figures. Therefore, the most straightforward strategy for classifying these polytopes
is as follows:
\begin{itemize}
\item Identify the projective and spherical rank $3$ polytopes.
\item For each pair of spherical or projective $\CK$ and $\CL$,
not both sperical, identify the universal $\{\CL,\CK\}$, if it exists.
\item Enumerate the quotients of each $\{\CL,\CK\}$, and identify those which are section regular.
\end{itemize}

Note for the second step, if $\CL$ is of type $\{p,q\}$ and $\CK$ of type $\{q',r\}$, a simple
necessary condition that $\{\CL,\CK\}$ exist is that $q$ must equal $q'$. Then $\{\CL,\CK\}$ 
will be of type $\{p,q,r\}$. As we shall see, this necessary condition is not sufficient.
The case $q=2$ cannot yield locally projective polytopes, since the polytopes
of type $\{p,2\}$ are ``dihedra'', with p-gonal faces, and in particular
are spherical (along with their duals). There do exist locally projective
polytopes with $p=2$ or $r=2$. The only cases for $r=2$ are ``ditopes''
$\{p,q,2\}$ with projective facets of type $\{p,q\}$. Their
enumeration is easy, and since the projective $\{p,q\}$ have no
proper quotients, neither do the ditopes. If $p=2$ we obtain the duals
of the ditopes. For most of the rest of this article, we shall ignore these
degenerate cases and assume $p,q,r\geq3$.

The spherical rank $3$ polytopes have been known since the time of the ancient Greeks. They are the
tetrahedron $\{3,3\}$, the cube $\{4,3\}$, the octahedron or cross $\{3,4\}$, the dodecahedron
$\{5,3\}$ and the icosahedron $\{3,5\}$.
Except for the tetrahedron, the groups of each of these
has a central inversion $\omega$. Taking the quotient of these polytopes by the group $\scl\omega$
yields four projective polytopes: the hemicube $\{4,3\}_3=\{4,3\}/2$, its dual the hemicross,
the hemidodecahedron $\{5,3\}_5=\{5,3\}/2$, and its dual the hemi-icosahedron. As noted in \cite{har3},
these are the only proper quotients of the platonic solids, and therefore, the only
rank $3$ section regular projective polytopes.

Ignoring duality leads to 22 cases that need to be analysed. These are listed in Table~\ref{cases}.

\begin{table}[ht]
\begin{center}
\begin{tabular}{|c|l|l|l||c|l|l|l|}
\hline
\#& Facets      & V.Figures   & Type         &\#& Facets      & V.Figures   & Type        \\
 \hline
1 & $\{3,3\}$   & $\{3,4\}_3$ & $\{3,3,4\}$  &2 & $\{3,3\}$   & $\{3,5\}_5$ & $\{3,3,5\}$ \\
3 & $\{3,4\}$   & $\{4,3\}_3$ & $\{3,4,3\}$  &4 & $\{3,4\}_3$ & $\{4,3\}_3$ & $\{3,4,3\}$ \\
5 & $\{3,4\}_3$ & $\{4,3\}$   & $\{3,4,3\}$  &6 & $\{3,5\}$   & $\{5,3\}_5$ & $\{3,5,3\}$ \\
7 & $\{3,5\}_5$ & $\{5,3\}_5$ & $\{3,5,3\}$  &8 & $\{3,5\}_5$ & $\{5,3\}$   & $\{3,5,3\}$ \\
9 & $\{4,3\}_3$ & $\{3,3\}$   & $\{4,3,3\}$  &10& $\{4,3\}$   & $\{3,4\}_3$ & $\{4,3,4\}$ \\
11& $\{4,3\}_3$ & $\{3,4\}_3$ & $\{4,3,4\}$  &12& $\{4,3\}_3$ & $\{3,4\}$   & $\{4,3,4\}$ \\
13& $\{4,3\}$   & $\{3,5\}_5$ & $\{4,3,5\}$  &14& $\{4,3\}_3$ & $\{3,5\}_5$ & $\{4,3,5\}$ \\
15& $\{4,3\}_3$ & $\{3,5\}$   & $\{4,3,5\}$  &16& $\{5,3\}_5$ & $\{3,3\}$   & $\{5,3,3\}$ \\
17& $\{5,3\}$   & $\{3,4\}_3$ & $\{5,3,4\}$  &18& $\{5,3\}_5$ & $\{3,4\}_3$ & $\{5,3,4\}$ \\
19& $\{5,3\}_5$ & $\{3,4\}$   & $\{5,3,4\}$  &20& $\{5,3\}$   & $\{3,5\}_5$ & $\{5,3,5\}$ \\
21& $\{5,3\}_5$ & $\{3,5\}_5$ & $\{5,3,5\}$  &22& $\{5,3\}_5$ & $\{3,5\}$   & $\{5,3,5\}$ \\
 \hline
 \end{tabular}
 \end{center}
 \caption{The Scope of the Classification Problem.}
 \label{cases}
 \end{table}

These cases have all been completely analysed, and the results reported in the literature, 
with the exception of cases 10 to 12,
for which a small amount of additional work remains to be done.

\section{The Unfinished Cases.}
\label{uncases}
\seq

In this short section, we consider cases 10 and 11. Case 12 is dual to case 10, 
and therefore does not warrant separate
consideration. As mentioned earlier, polytopes of these types exist (see \cite{amalg}). Furthermore,
the only polytopes with the given facets and vertex figures are in fact the universal polytopes,
as was shown in \cite{har5}. However, a full analysis of the quotients of these polytopes has not 
been done. It will be useful at this point to review the classification of the quotients of the cube
(and therefore the cross), given in \cite{har3}. Let the group of the cube be 
$W'=\scl{s_0,s_1,s_2}$, and let $x=s_0$, $y=s_1xs_1$ and $z=s_2ys_2$. The semisparse
subgroups of $W'$ are the trivial subgroup (leading to the cube itself), the group
$\scl{xy}$ and its conjugates $\scl{yz}$ and $\scl{xz}$ (leading to the digonal prism),
the group $\scl{xyz}$ (leading to the hemicube), and the group $\scl{xy,yz}$ (leading 
to the polytope $\{2,3\}$).

In particular, note that the $3$-hemicube (and therefore the $3$-hemicross) has no proper quotients
(this is not true in higher ranks). Since the vertex figures of $\{\{4,3\},\{3,4\}/2\}$ have no
proper quotients, Theorem 2.7 of \cite{har6} may be applied to case 10, so that semisparse subgroups
of $W=\scl{s_0,s_1,s_2,s_3}=[\{4,3\},\{3,4\}/2]$ are characterised by the property that all their
conjugates intersect $\scl{s_0,s_1,s_2}\scl{s_1,s_2,s_3}$ in a semisparse subgroup of $W'$.
It is simple enough (although a little unsatisfying) to do a computer
search of the subgroups of $W$ to find all those which satisfy the required
property. This leads to the following theorem and corollary.

\bgth
\label{quo434}
\thm The universal polytope $\CP=\{\{4,3\},\{3,4\}_3\}$ has four quotients. These are $\CP$ itself,
$\{\{4,3\}_3,\{3,4\}_3\}$, $\{\{2,3\},\{3,4\}_3\}$, and a nonregular polytope whose facets 
are all digons.
\ndth

Note that the four quotients of $\CP$ correspond very naturally to the four quotients of the cube.
Note also the following Corollary.

\bgth
\label{quo434b}
\cor The polytope $\CQ=\{\{4,3\}_3,\{3,4\}_3\}$ has no proper quotients.
\ndth

The author notes that it would be a worthwhile endeavour to characterize the quotients of
$\{\{4,3^{n-3}\},\{3^{n-3},4\}/2\}$ for general $n\geq4$, but feels that such a characterization
falls beyond the scope of this article.

%First, consider the rank $n$ polytope $\CP=\{\{4,3^{n-3}\},\{3^{n-3},4\}/2\}$, the universal 
%polytope with cubes for facets and hemicrosses for vertex figures. 
%As noted in \cite{har5}, this polytope is flat, and has 
%$2^{n-2}$ vertices and the same number of facets. Let $N$ be a semisparse subgroup of its group
%$W$, let $H_0$ be the group of the vertex figures, $H_3$ the group of the facets of $\CP$ and denote
%the group of $\CP$ by $W$. The polytope $\CP$ may be constructed as $2^\CH$, where $\CH$ is the hemicross,
%using the twisting operation developed by McMullen and Schulte (see Chapter 9 of \cite{msbook}).

\section{The Classification.}
\label{class}
\seq

This final section of the article will go through the cases of Table~\ref{cases} one by one, giving 
a description of the universal polytope and its quotients, with references to the literature.

\begin{itemize}
\item Cases 1 to 5: The results of \cite{har3}, where polytopes of ``finite
type'' were studied (that is, polytopes with the same Schl\"afli symbols
as the classical spherical polytopes) show that there are no locally
projective polytopes of types $\{3,3,4\}$, $\{3,3,5\}$ or $\{3,4,3\}$. This
likewise eliminates cases 9 and 16, as is noted below.
\item Case 6: There is no polytope
with icosahedral facets and hemidodecahedral vertex figures. See the dual Case 8 for an
explanation.
\item Case 7 is Coxeter and Gr\"unbaum's 11-cell (see \cite{cox11} or \cite{gru11}).
The polytope has a group of order 660, isomorphic to the projective special linear
group $L_2(11)$. It was noted in \cite{har4} that this polytope has no proper
quotients.
\item Case 8: There is no polytope with hemiicosahedral facets and dodecahedral vertex figures. If
such a polytope existed, it could be constructed (in principle) by arranging hemiicosahedra face to 
face, with three around each edge, either until the polytope closed up, or ad infinitum. However,
an attempt to perform this construction results in the 11-cell of Case 7, as noted in \cite{gru11}.
\item Case 9 is dual to case 1, thus yields no polytopes.
\item Case 10: The universal polytope $\{\{4,3\},\{3,4\}_3\}$ appeared in
\cite{amalg}. A simple proof that the universal polytope is the only polytope
of type $\{\{4,3\},\{3,4\}_3\}$ appeared in \cite{har5}. The polytope may be 
constructed from its vertex figures via the ``twisting'' operation developed by
Schulte and McMullen and described in Chapter 9 of \cite{msbook}. Here, 
$\{\{4,3\},\{3,4\}_3\}$ is $2^\CH$, where $\CH=\{3,4\}_3$ is the hemicross.
Its group is the semidirect
product $2^3\rtimes S_4$, with $S_4$ (the group of $\CH$) acting by conjugation on the 
generators of $2^3$ in the same way that it acts on the vertices of the hemicross.
The quotients of this polytope were described in Theorem~\ref{quo434}.
\item Case 11: The universal polytope $\{\{4,3\}_3,\{3,4\}_3\}$ is a quotient of
$\{\{4,3\},\{3,4\}_3\}$ and has itself no proper quotients. This was shown in 
Theorem~\ref{quo434b}, but could easily have been noted from the results of 
\cite{har5}. Its group is the quotient of $2^3\rtimes S_4$ by the
(normal) subgroup generated by the product of the three generators
of $2^3$.
\item Case 12 is dual to case 10.
\item Case 13 was studied in depth in \cite{har4}, where it was shown that the universal
$\{\{4,3\},\{3,5\}_5\}$ exists and is finite. Its existance and finiteness are
also noted
in Section 8E of \cite{msbook}. It appears there as $2^\CI$, where $\CI$ is the hemi-icosahedron
(whose group is $A_5$). Its group is therefore $2^6\rtimes A_5$, with $A_5$ acting by conjugation
on the generators of $2^6$ in the same way that it acts on the vertices of $\CI$. 
The universal polytope therefore has 80 cubes as facets, and 64 vertices. In
\cite{har4} it was shown, via a computer search, that the polytope has 70 quotients. Of these,
three are regular, these being $\{\{2,3\},\{3,5\}_5\}$, the universal polytope, and a quotient
of the universal polytope by a group of order $2$. Besides the latter two, there are another 
nine which are nonregular polytopes of type $\{\{4,3\},\{3,5\}_5\}$, and besides these
another eight of type $\{4,3,5\}$. The latter eight each have, as facets, some cubes and 
some hemicubes, hence they are not section regular.
For more details on the quotients of this polytope, the reader is referred
to \cite{har4}.
\item Cases 14 and 15: It was shown in Theorem 3.6 of \cite{har5} that a polytope with
hemicubes for facets cannot be of type $\{4,3,...,3,p\}$ for odd $p$. Thus cases 14 and 15 do not
yield examples of locally projective polytopes. 
\item Case 16 is dual to case 2.
\item Cases 17, 18 and 19 are dual to cases 15, 14 and 13 respectively.
\item Case 20 was examined in \cite{harlee}. The polytope $\{\{5,3\},\{3,5\}_5\}$ exists and is
finite, with group $J_1\times L_2(19)$ of order 600415200. It therefore has 10006920 vertices and
half that many facets. Here, $L_2(19)$ is a projective special linear group of order 3420, and
$J_1$ is the first Janko group, of order 175560. The quotients of this polytope are 145 in
number, including as regular quotients Coxeter's 57-cell $\{\{5,3\}_5,\{3,5\}_5\}$ (see below),
and a polytope of type $\{\{5,3\},\{3,5\}_5\}$ with group $J_1$. There are also
67 nonregular quotients of this type. The remaining 75 proper quotients
(all of type $\{5,3,5\}$) each have, as facets, some
dodecahedra and some hemidodecahedra. For more details, see \cite{harlee}.
\item Case 21 was discovered by Coxeter in \cite{cox57}. The universal polytope has 
57 facets, and its group is $L_2(19)$ of order 3420. It was shown in \cite{har4} that it has
no proper quotients. It is itself a quotient of the universal polytope of case 20.
\item Case 22 is dual to case 20, and therefore does not warrant separate discussion.
\end{itemize}

To summarize, there are nine nondegenerate 
universal locally projective regular polytopes of rank $4$, namely 
$\{\{3,5\}_5,\{5,3\}_5\}$, $\{\{4,3\},\{3,4\}_3\}$ and its dual,
$\{\{4,3\}_3,\{3,4\}_3\}$, $\{\{4,3\},\{3,5\}_5\}$ and its dual, 
$\{\{5,3\}_5,\{3,5\}_5\}$ and $\{\{5,3\},\{3,5\}_5\}$ and its dual.
There are 13 regular polytopes of these types, since each of those with hemi-icosahedral
vertex figures (or hemidodecahedral facets) has a proper regular quotient of the same type.

These universal polytopes have a total of $437$ quotients, 17 regular, and 169 section regular.
The extra four regular quotients are the quotients $\{\{2,3\},\CK\}$ of locally projective\
$\{\{4,3\},\CK\}$ and their duals. These are also counted amongst the section regular quotients.
The 152 nonregular section regular quotients (or their duals) all have hemi-icosahedral vertex figures
and spherical (cubic or dodecahderal) facets. Besides these, there are a further
four degenerate locally projective polytopes, namely $\{\{2,4\},\{4,3\}_3\}$ and
$\{\{2,5\},\{5,3\}_5\}$ and their duals. These four are not quotients of any
nondegenerate locally projective polytopes.

The most remarkable feature of the classification is the following result.

\bgth
\label{finite}
\thm All locally projective rank 4 polytopes are finite.
\ndth

The author would like to acknowledge and thank Jesus Christ, through whom
all things were made, for the encouragement and motivation to 
complete this article.


\begin{thebibliography}{99}

%\bibitem{magma} W. Bosma, J. Cannon, C. Playoust, ``The Magma Algebra System I: the User Language'',
%J. Symbolic Comput. {\bf 24}, 235--265 (1997).

%\bibitem{bou} N. Bourbaki, ``Groupes et Alg\`ebres de Lie'', Chapitres IV--VI (Hermann, 1968).

\bibitem{cox} H. S. M. Coxeter, ``Regular Polytopes'', (Methuen and Co., 1968).

\bibitem{cox57} H. S. M. Coxeter, ``Ten Toroids and Fifty-Seven Hemi-Dodecahedra''
Geom. Dedicata {\bf 13}, 87--99 (1982).

\bibitem{cox11} H. S. M. Coxeter, ``A Symmetrical Arrangement of
Eleven Hemi-icosahedra'', Ann. Disc. Math. {\bf 20}, 103--114 (1984).

%\bibitem{gap4} The GAP Group, ``GAP~-- Groups, Algorithms, and Programming,
%Version 4.3'', http://www.gap-system.org (2002).

%\bibitem{gev} G. G\'evay, ``On Perfect 4-Polytopes'', Beit. Alg. Geom. {\bf 43},
%243--259 (2002) 0138-4821/93.

\bibitem{gru11} B. Gr\"unbaum, ``Regularity of Graphs, Complexes and Designs'',
In: ``Probl\`emes Combinatoires et Th\'eorie des Graphes, Colloquium International
CNRS, Orsay'' {\bf 260}, 191--197 (1977).

%\bibitem{har} M. I. Hartley, ``Combinatorially Regular Euler Polytopes'',
%PhD Dissertation (University of Western Australia, 1996).

%\bibitem{har1} M. I. Hartley, ``All Polytopes are Quotients, and
%Isomorphic Polytopes are Quotients by Conjugate Subgroups'',
%Discrete Comput. Geom. {\bf 21}, 289--298 (1999).

%\bibitem{har2} M. I. Hartley, ``More on Quotient Polytopes'',
%Aequationes Math. {\bf 57}, 108--120 (1999), 0001-9054/99/0101080-13.

\bibitem{har3} M. I. Hartley, ``Polytopes of Finite Type'',
Discrete Math. {\bf 218}, 97--108 (2000) PII: S0012-365X(99)00339-8.

\bibitem{har4} M. I. Hartley, ``Quotients of Some Finite Universal
Locally Projective Polytopes'', Discrete Comput. Geom. {bf 29}, 435--443 (2003) 
DOI: 10.1007/s00454-002-2852-y.

\bibitem{har5} M. I. Hartley, ``Locally Projective Polytopes of Type $\{4,3,\dots,3,p\}$'',
Submitted to: Journal of Algebra.

\bibitem{har6} M. I. Hartley, ``Simpler Tests for Semisparse Subgroups'',
Submitted to: Geometry and Topology.

\bibitem{harlee} M. I. Hartley, D. Leemans, ``Quotients of a Locally Projective 
Polytope of Type $\{5,3,5\}$. To Appear: Math. Zeit.

%\bibitem{535aux} M. I. Hartley, D. Leemans, ``Quotients of a Universal Locally
%Projective Polytope of type $\{5,3,5\}$ (Auxiliary Information)'',
%http://www.angelfire.com/mt/ofolives/535.html (2003).

%\bibitem{hum} J. E. Humphreys, ``Reflection Groups and Coxeter Groups'',
%Cambridge Studies in Advanced Mathematics (Cambridge University Press, 1990).

%\bibitem{janko} Z. Janko, ``A New Finite Simple Group withe Abelian
%Sylow 2-subgroups, and its Characterisation'', J. Algebra {\bf 3}, 147--186 (1966).

%\bibitem{liv} D. Livingstone, ``On a Permutation Representation of the Janko Group'',
%J. Algebra {\bf 6}, 43--55 (1967).

\bibitem{mslocproj} P. McMullen, ``Locally Projective Regular Polytopes'',
J. Comb. Theory A {\bf 65}, 1--10 (1994).

%\bibitem{mcq} P. McMullen, E. Schulte, ``Quotients of Polytopes and C-groups'',
%Discrete Comput. Geom. {\bf 11}, 453--464 (1994).

%\bibitem{msflat} P. McMullen, E. Schulte, ``Flat Regular Polytopes'',
%Annals of Combinatorics {\bf 1} 261--278 (1997).

\bibitem{msbook} P. McMullen, E. Schulte, ``Abstract Regular Polytopes''
(Cambridge University Press, 2002).

\bibitem{amalg} E. Schulte, ``Amalgamation Of Regular Incidence-Polytopes'',
Proc. London Math. Soc. {\bf 56}, 303--328 (1988).

\bibitem{mstwis} E. Schulte, ``Classification of Locally Toroidal Regular
Polytopes'', in: T. Bisztriczky et al. (Eds.) ``Polytopes: Abstract,
Convex and Computational'', 125--154 (Kluwer, 1994).

%\bibitem{weiss} A. I. Weiss, ``Incidence Polytopes of Type $\{6,3,3\}$'',
%Geom. Dedicata {\bf 20} 147--155 (1986).





\end{thebibliography}
\end{document}